\documentstyle[amscd]{amsart} 
\DeclareMathSymbol{\epi}  {\mathrel}{AMSa}{"10}   

\newtheorem{thm}{Theorem}[section]
\newtheorem{prop}[thm]{Proposition}
\newtheorem{lem}[thm]{Lemma}
\newtheorem{cor}[thm]{Corollary}

\theoremstyle{definition}
\newtheorem{defi}[thm]{Definition}
\newtheorem{rem}[thm]{Remark}

\newtheorem{exa}[thm]{Example}

\theoremstyle{remark}

\font \Bbb=msbm10 scaled \magstep1

\def \N{\hbox{\Bbb N}}

\def \>{\rangle}
\def \<{\langle}
\def \Xm{X_1,\dots,X_m}

\def \op{\operatorname}
\def \max{\op{max}}
\def \min{\op{min}}
\def \cd{\op{crit.deg.}}
\def \dim{\op{dim}}
\def \soc{\op{soc}}
  
\begin{document}

\title{Artinian Algebras and differential forms}

\author{Guillermo Corti\~nas}
\author{Fabiana Krongold}
\address{Depto.de Matem\'atica\\
	 Ciudad Universitaria Pab. I\\ 
	 (1428) Buenos Aires\\ Argentina}
\email{gcorti@@mate.dm.uba.ar, fkrong@@mate.dm.uba.ar}
\maketitle

%

\setcounter{section}{-1}
\begin{section}{Introduction}

This article concerns commutative algebras over a field $k$ of 
characteristic zero which are finite dimensional as vectorspaces, and
particularly those of such algebras which are graded. 
Here the term graded is applied to non-negatively graded 
algebras $A$ with $A_0$ reduced and finite dimensional. Thus the trivial 
grading $A=A_0$ is only allowed if $A$ is a product of finite field 
extensions of $k$. It was conjectured in [2] that for all finite dimensional
algebras $A$ which are not principal ideal algebras (i.e. have at least one
nonprincipal ideal), 
the following submodule of the K\"ahler 
differentials is nonzero: 
$$\tau(A)=\bigcap{\ker(\Omega_A @>>>\Omega_B)}$$
Here the intersection is taken over all principal ideal algebras $B$ 
 and all homomorphisms 
$A @>>>B$. In this paper we prove that the conjecture holds for both 
Gorenstein graded and standard graded algebras. The conjecture,
called Artinian Berger Conjecture (ABC) in [2], was introduced in
connection with the classical Berger Conjecture (BC) of [1]. In fact it
is proven [2, Theorem 0.1] that ABC $\Longrightarrow$ BC. 
 (Both the BC and the result of [2] are recalled in 1.2 below.)
However we shall show here that the implication 
is not as straightforward as one could expect. For instance our theorem that ABC
holds for standard graded algebras (2.1) implies BC for some graded and some 
ungradable algebras of dimension one (2.2). Our theorem that
0-dimensional Gorenstein graded algebras satisfy ABC (3.2) proves no case of 
the BC
(3.3); in particular it does not prove the latter for 1-dimensional Gorenstein
graded algebras.

\medskip
The ABC is connected to other interesting questions such as when is an 
artinian algebra embeddable in a principal ideal algebra. For instance we 
show that nonprincipal Gorenstein algebras are not embeddable (\ref{gor})
and we classify standard graded algebras in terms of their degree of 
embeddability (2.5).

\medskip
The rest of this paper is organized as follows. Some basic properties from
[2] are recalled in section 1. Sections 2 and 3 are devoted respectively
to graded standard and Gorenstein algebras. Theorem 2.1 and Proposition 2.4 
were obtained jointly with P. Solern\'o. Theorem 3.1 was proven 
independently --and almost simultaneously-- by S. Geller and C. Weibel. 

\end{section}

\bigskip
\begin{section}{Preliminaries}
The two main results of this paper are stated for local algebras over an
arbitrary field of characteristic zero. The following lemma allows us to 
prove them in the
algebraically closed case only. By part ii) of the lemma bellow we can 
rephrase the main results of this paper  (\ref{abc} and \ref{gor}) into one;
namely that ABC holds for products of Gorenstein graded and of standard 
graded algebras. 

\begin{lem}
Let $A,A_1$ and $A_2$ be $k$-algebras; write ${\overline k}$ for the algebraic
closure of $k$. Then:

\item 
i) If ABC holds for $A\otimes_k{\overline k}$ then it holds for $A$.
\item
ii) If ABC holds for $A_1$ and $A_2$ then it holds for the product 
$A_1\times A_2$.

\end{lem}

\begin{pf} 
It is immediate from [2, 1.1.1 and 2.0]
\end{pf}

{\bf Artinian Berger Conjecture - Berger Conjecture}:
The main result of [2] states that if ABC holds for every 0-Krull 
dimensional $k$-algebra, then BC holds for every 1-Krull dimensional reduced
algebra. Actually the proof given in [2] proves the following more precise
statement:

\begin{lem}\label{pb}
Suppose a pull back diagram:
\begin{equation} \label{pbd}
\begin{CD}
R   @>>> S\\
@VVV     @VVV\\
A   @>>> B
\end{CD}
\end{equation}
is given where $A \hookrightarrow B$ is an inclusion of 0-Krull dimensional 
algebras, the vertical maps are surjective, $A$ is not a principal 
ideal algebra, $S$ is reduced normal and Krull-dim$S=1$. Then $R$ is
reduced singular, Krull-dim$R=1$ and $S$ is its normalization.
If $\Omega_A\to \Omega_B$ is not injective then $\Omega_R@>>>\Omega_S$ is
not injective either, i.e. BC holds for  $R$.
\end{lem}

\begin{pf}
That $R$ satisfies the hypothesis of BC is not hard to see. If 
$\Omega_A\to \Omega_B$ is not injective then $\Omega_R\to \Omega_S$
cannot be injective, by [2, 1.1.1]
\end{pf}

\begin{rem}
Note that if $A$ occurs in a diagram such as (\ref{pbd}) then it is necessarily
embeddable. It follows from the proof of [2, 2.3] that the canonical map:
\begin{equation}\label{dr}
H_{dR}^0(A)@>>>H_{dR}^0(A_{red})
\end{equation}
is an isomorphism if $A$ is embeddable. 
 On the other hand if $A$ is not embeddable and (\ref{dr}) is an 
isomorphism, then it satisfies ABC.
To see this, first use 1.1 ii) to restrict to the local case. Next let
$(A,M)$ be local and unembeddable, and let $0\neq m$ such that any map
$A\to B$
into a principal ideal algebra maps it to zero. Then its differential 
$dm\in \tau(A)$ is a non-zero element by (\ref{dr}). It should be noted
that the map (\ref{dr}) need not be an isomorphism in general, even 
under our standing assumption that $char(k)=0$, as shown by the example
below.  What our assumption that $char(k)=0$ does imply is that  
(\ref{dr}) holds when $A$ is graded --in the sense of this paper; see 
\S 0 . Indeed in this case the derivation $D(a)=deg(a)a$ is injective on
$A_+$, whence so must be $d:A@>>>\Omega_A$. 

\end{rem}

\begin{exa}
Consider the ideal  $I:=<X^3Y, X^5,XY^3+2X^3,3X^2Y^2+5Y^4>\subset k[X,Y]$;
note $I$ contains both partial derivatives of $F=X^4+X^2Y^3+Y^5$. Hence
the image $f$ of $F$ in $A=k[X,Y]/I$ is in 
$K=\ker(H_{dR}^0(A)@>>>H_{dR}^0(A_{red})=k)$. One checks
that $\{ 1,x,y,x^2,xy,y^2,x^3,x^2y,xy^2,y^3,x^4,x^2y^2\}$ is a basis
of the vectorspace $A$ and that $f=\frac{x^4}{5}$, whence it is a nonzero
element of $K$. As per the remark 
above, this implies that $A$ is unembeddable. To see that ABC actually 
holds for $A$, proceed as follows. First one checks that the basis 
element $w=x^2y^2$ is mapped to zero by any map into a principal ideal 
algebra. Hence $dw\in\tau(A)$. Further, $dw\ne 0$, since $w$ goes to an
element of positive degree in the graded algebra $A/<x^3>$.

\end{exa}

In order to prove the main results of this paper we shall use the following 
artinian version of valuation theory.

\medskip
{\bf Truncated Valuations:} Let $\alpha :A@>>>B=k[t]/\<t^{N+1}\>$ be a map. Consider 
the map  $\nu=\nu_{\alpha}:A@>>>\{0,1,\dots,N,\infty\}$ given by
$\nu(a)=e$ when $\alpha(a)=t^eu$ for some invertible element 
$u\in B$, and $\nu(a)=\infty$ if $\alpha(a)=0$.
We call $\nu$ the {\em truncated valuation \/} associated to $\alpha$.
Note that truncated valuations are not valuations in the classical sense
of [5, Ch. VI,\S 8, page 32]. However our truncated valuations do retain 
some of the properties
of classical valuations. For instance we have
$\nu(xy)=\nu(x) + \nu(y)$ and $\nu(x+y) \le \min\{\nu(x),\nu(y)\}$. Thus
$\nu$ is a morphism from the  
non-cancellative multiplicative monoid $(A,\cdot)$ to the additive monoid
$\{0,1,\dots,N,\infty\}$, where $n+\infty=\infty+n=\infty$. The order in 
$\{0,1,\dots,N,\infty\}$ is the obvious one. Another trivial, but useful 
property of truncated valuations is the following. Let 
$a_1,\dots, a_n\in A$ be a $k$-linearly independent set,
$S=\sum ka_i$, and $d=n-\dim\ker\alpha\cap S$. Assume $d\ne 0$, i.e.
assume $S\not\subset\ker\alpha$. 
Then there exist elements $x_1,\dots ,x_n\in A$
such that: 
\begin{equation}\label{vx}
\nu(x_1)<\nu(x_2)\dots<\nu(x_d),{\ \ }\nu(x_i)=\infty {\ \ }(i>d)
\end{equation}
and such that still $\sum kx_i=S$. This property follows from the fact
that the $n\times N$ matrix having as rows the coefficients of 
$\alpha(a_1),\dots,\alpha(a_n)$ can be made into an upper triangular
matrix through row operations.

\end{section}
\bigskip 
\begin{section}{Standard Graded Algebras}

Let  $A= A_0\oplus A_1 \oplus \dots \oplus A_n$\ be a finite dimensional
graded algebra. We say that $A$ is {\em standard\/} if the ideal
$M=A_1 \oplus \dots \oplus A_n$ is maximal and is generated by degree one 
elements. Thus $(A,M)$ is a local ring with residue field $A_0$. Hence
any standard graded algebra may be written as $A=F[\Xm]/I$ where
$F$ is a finite field extension of the ground field $k$, and $I$ is 
homogeneous for the standard grading of the polynomial ring.\\
The main theorem of this section is the following:

\begin{thm}\label{abc}
Standard graded algebras satisfy ABC.
\end{thm}

The next example shows how ungraded cases of BC can be deduced from
the  theorem above.

\begin{exa}
Consider the smooth curve of genus one $S=k[X,Y]/\<X(X^2-1)+Y^3\>$ and the
standard graded algebra $A=k[X,Y]/\<X,Y\>^4$. Form the pull-back:
$$
\begin{CD}
R   @>>> S\\
@VVV     @VVV\\
A   @>>> B
\end{CD}
$$
Here $B=k[t]/\<t^{16}\>=S/\<X,Y\>^{16}$ and $A@>>>B$ is the inclusion
$X\mapsto t^4, Y\mapsto t^5$. It is not hard to see --using [4]-- that 
$S$ is not gradable whence $R$ cannot be gradable either.  Thus by 
lemma~\ref{pb} and the theorem above, BC holds for the ungradable ring $R$.
\end{exa}

\bigskip
The proof of theorem \ref{abc} will take us the rest of this section. Since 
standard grading is preserved by change of ground field, the proof of the 
theorem immediately reduces to the case when $k$ is algebraically closed. 
For the rest of this section we shall assume $k$ is algebraically closed.

\smallskip
It will be essential for our proof to measure the degree of embeddability of 
standard graded algebras in truncated polynomial rings. The critical degree
defined in 2.4 below will be a useful invariant for this purpose.

\begin{lem} \label{dl}
Let $(A,M)$ be a standard graded algebra. If $A$ is not a principal ideal 
algebra then the following is a finite nonempty set:
$${\cal E}=\{i\in \N / \ \exists\ B \mbox{\ truncated polynomial 
algebra and\ } \alpha :A \to B \ \mbox{\ s.t.\ }\dim\alpha(A_i) \ge 2\}$$

\end{lem}

\begin{pf}
Let  $n\ge 1$ be such that $M^{n+1}=0$.
Then ${{\cal E}\subset\{1,\dots,n\}}$. Further, $1\in\cal E$ because there 
exists a map $A/M^2\to k[t]/\<t^4\>$ sending $M/M^2$ onto $\<t^2,t^3\>$ since 
$\dim M/M^2\ge 2$.
\end{pf}

\begin{defi}
In the situation of the lemma above, we define the {\em critical degree\/} of $A$ as
the integer $\cd(A)=\max({\cal E})$. 
\end{defi}

\begin{prop}\label{epi}
Let $A$ be a standard graded algebra. Let $r$ be its critical degree.
Then there exists a surjective homomorphism:
$$A \epi Q(r)=k[X,Y]/\<X^{r+1},X^rY,Y^2\>$$
\end{prop}

\begin{pf}
Let $\alpha:A\to B$ be a map with values in a truncated polynomial algebra 
such that dim$\alpha(A_r)\ge 2$. Let $\{x,y,z_3,\dots,z_m\}$ be a basis of 
$A_1$ such that ${e=\nu(x)<f=\nu(y)}$ are the least two valuations in $A_1$, 
and such that $\nu(z_i)>f,(i=3,\dots,n)$. Consider the ideal 
${I=\<x^{r+1},x^ry,y^2,z_3,\dots,z_m\>\!\subset\!A}$. 
We have a surjective map $\phi:Q(r)\epi A/I, X\!\mapsto\!x, 
Y\!\mapsto\!y$.
A valuation argument shows that for each $1\!\le\!i\!\le\!r$ the elements 
${x^i,x^{i-1}y}$ are linearly independent in $(A/I)_i$. 
It follows that $\phi$ is an isomorphism.
\end{pf}

\begin{rem}
It is not hard to show --using valuations-- that $\cd Q(r)=[{r\over 2}]+1$
(here $[x]$ is the integer part of $x$).
Hence the fact that an epimorphism ${A\epi Q(r)}$ exists does not imply any 
relation between the critical degree of $A$ and $r$. 
In other words, $\cd(A)\le \max\{r:\exists A\to Q(r)\}$, but the
inequality may be strict. 
\end{rem}

\medskip
\begin{pf*}{Proof of Theorem~\ref{abc}}
Assume $A$ is not a principal ideal algebra
and let $r$ be its critical degree. We shall show that for all $x,y \in A_1$
the element:
\begin{equation}
\omega(x,y)=x^{r-1}(xdy-ydx)\in \tau(A)  \label{eq}
\end{equation}
On the other hand it is not hard to show that if $x,y$ are as in the proof 
of \ref{epi} above then composite
$A\epi A/I@>{\phi^{-1}}>> Q(r)$  maps $\omega$ to a non-zero element.
Thus it suffices to prove (\ref{eq}).
Let $\beta:A\to k[t]/\<t^{N+1}\>$ be any map. We have to show that 
$\beta(\omega)=0$. Write $e$ and $f$ for the valuations of $x$ and $y$,
and ${e'=\min\{\nu(a_1) : a_1\in A_1\}}$. 
We have that ${\nu(\omega)\ge re+f}$. If ${re+f= \infty}$ we are done. 
Otherwise ${0 \ne \beta(x^ry)\in \beta(A_{r+1})}$. Then 
${\beta(A_{r+1})=k\beta(x^ry)}$ and therefore ${re+f=(r+1)e'}$, whence 
${e=f=e'}$. Choose a basis ${\{x,z_2,\dots,z_m\}}$ of $A_1$ with 
${\nu(z_i)>\nu(x)}$. Then $y=\lambda x+p$ where $\lambda \in k$ and 
$p=\sum \lambda_iz_i,\lambda_i \in k$. As $\omega(x,-)$ is linear, 
$\omega(x,y)=\lambda\omega(x,x)+\omega(x,p)$.
But now ${\nu(p) > \nu(x)}$, whence 
$\beta(\omega(x,y)=\beta(\omega(x,p))=0$ by the argument above.
\end{pf*}
\end{section}
\bigskip
\begin{section}{Gorenstein Algebras}
Recall a finite dimensional local algebra $(A,M)$ is called {\em Gorenstein\/} 
(or self-injective) if for the socle we have $\dim_k\soc(A)=1$.

\begin{thm} \label{gor}
Let $A$ be a finite dimensional Gorenstein local
algebra over a field $k$. Assume $A$ is not a principal ideal algebra.
Then every homomorphism from $A$ into a principal ideal algebra maps the
socle of $A$ to zero. In particular $A$ is unembeddable.
\end{thm}

\begin{pf}
We may assume $k$ is algebraically closed (by 1.1). Let $M$ 
be the maximal ideal of $A$ and let $n$ be a positive integer such that 
$M^n\ne 0$
and $M^{n+1} =0$. As $A$ is Gorenstein, $M^n=\soc A$. If 
a homomorphism to a product of truncated polynomial algebras which is
not zero on $M^n$ exists, 
then there must also exist a
homomorphism $f:A@>>> k[t]/\<t^{N+1}\>$ such that $f(M^n)\ne 0$.
Because $\dim M^n=1$, $f$ must be injective on $M^n$. 
Let $\{x_1,\dots , x_m\} \subset M$ a minimal generating set of $M$ 
with ascending valuations as in (\ref{vx}) above.
Then $f(x_1^n) \ne 0$ and 
any other monomial of degree $n$ in the $x_i$ is in $\ker f$. 
As $f$ is injective 
on $M^n$  we have $M^n=kx_1^n$, and any other monomial in the $x_i$
is zero. We have proved the following statement for $j=0$:
if $n-j \ge 1$ then $M^{n-j} = kx_1^{n-j}+\dots +kx_1^n$. 
Assume by 
induction that this assertion is true for $j$. If $n-j=1$ we are done. 
Otherwise $n-(j+1) \ge 1$. Let $\alpha \in M^{n-(j+1)}$; we must prove 
that $\alpha \in kx_1^{n-(j+1)} + \dots + kx_1^n$. Since 
$x_1\alpha \in M^{n-j}$, then
 $x_1\alpha = \sum_{i=0}^{n-(j+1)}\lambda_ix_1^{n-i}$, by induction.
Write $\beta=\alpha - \sum_{i=0}^{n-j}\lambda_ix_1^{n-i-1}$. We have 
$x_1\beta=0$. Hence $\infty=\nu(x_1\beta) = \nu(x_1) + \nu(\beta) \le
\nu(x_i) + \nu(\beta) =\infty$ for $1<i\le m$. We claim $x_i\beta=0$ in 
$A$. To see this note first that $x_i\beta = \sum_{l=0}^j 
\mu_{il}x_1^{n-l}$. Next, apply
$f$ to both sides of this identity to get that $\mu_{il}=0$ for all $l$. 
Hence
$x_i\beta =0$ as claimed. It follows that $\beta \in M^n$ whence $\beta = 
\lambda x_1^n$ for some $\lambda$. Then $\alpha \in kx_1^{n-j-1}+\dots 
+kx_1^n$ as we wanted to prove. 
\end{pf}

\begin{cor}
Let $A=A_0\oplus A_1\oplus A_2\oplus\dots$ be a finite dimensional
graded algebra. Assume $(A,A_+)$ is a local Gorenstein ring. Then 
$A$ satisfies ABC.
\end{cor}

\begin{pf}
Immediate from the theorem and from remark 1.3.
\end{pf}
\begin{rem}
The result above does not give any information on the one dimensional BC.
Indeed by theorem~\ref{gor} above, no zero dimensional Gorenstein algebra 
$A$ can occur in the pullback diagram (\ref{pbd}). Also note
that the hypothesis that $A$ be graded can be replaced by the hypothesis
that (\ref{dr}) be an isomorphism, or even that $d(soc(A))\ne 0$. The
latter
condition is exactly what is needed to make the proof above work, and
is not automatic. Indeed one checks that the algebra $A$ of example
1.4 has $soc(A)=kx^4$; it was shown already that $d(x^4)=0$, whence $A$ is 
an ungradable Gorenstein algebra for which the proof above does not work. 
On the other hand ABC does hold for $A$, as shown in 1.4 above.
\end{rem}
\end{section}
\smallskip
{\it Acknowledgements:} An earlier version of this 
paper wrongly stated that the map (\ref{dr}) is always an isomorphism in 
characteristic zero. We are especially indebted to 
S. Geller and C. Weibel for noting the error and for bringing  example
1.4 to our attention. The latter, due to Weibel, is a variation of the 
example appearing in page 241 of [3]. Thanks also to L. Avramov from 
whom we learned of Saito's paper [4].


\smallskip
\begin{thebibliography}{10}
\bibitem{}
R. Berger, {\em Differentialmoduln eindimensionaler lokaler Ringe,}
Math. Zeit. {\bf 81} (1963), 326-354.

\bibitem{}
G.Corti\~nas, S.Geller, C.Weibel, {\em The Artinian Berger Conjecture,} Math. Zeit., to appear. 

\bibitem{}
H. Grauert, H. Kerner, {\em Deformationen von Singularit\"aten komplexer
R\"aume,} Math. Annalen {\bf 153} (1964), 263-260.

\bibitem{}
K.Saito, {\em Quasihomogene isolierte Singularit\"aten von Hyperfl\"achen,}
Invent. Math. {\bf 14,} (1971) 123-142. 

\bibitem{}
O. Zariski, P. Samuel, {\em Commutative algebra, Volume II,} Grad.
Texts in Math. {\bf 29}, Springer Verlag 1960.

\end{thebibliography}
\end{document}